\documentclass[a4paper,11pt]{amsart}
\usepackage{graphicx}
\usepackage{mathptmx}
\usepackage{amsmath}
\usepackage{amssymb}
\usepackage{enumitem}
\usepackage{xcolor}
\usepackage{xparse}
\NewDocumentCommand{\eulerian}{omm}
 {%
  \genfrac<>{0pt}{}{#2}{#3}%
  \IfValueT{#1}{_{\!#1}}%
 }

 \newmuskip\pFqmuskip

\newcommand*\pFq[6][8]{%
  \begingroup 
  \pFqmuskip=#1mu\relax
  \mathchardef\normalcomma=\mathcode`,
  \mathcode`\,=\string"8000
  \begingroup\lccode`\~=`\,
  \lowercase{\endgroup\let~}\pFqcomma
  {}_{#2}F_{#3}{\left(\genfrac..{0pt}{}{#4}{#5}\bigg|#6\right)}%
  \endgroup
}
\newcommand{\pFqcomma}{{\normalcomma}\mskip\pFqmuskip}

\newtheorem{theorem}{Theorem}
\newtheorem{lemma}[theorem]{Lemma}
\newtheorem{corollary}[theorem]{Corollary}

\begin{document}

\title[A Study on degenerate Whitney numbers of the first and second kinds]{A Study on degenerate Whitney numbers of the first and second kinds of Dowling lattices}

\author{Taekyun  Kim}
\address{Department of Mathematics, Kwangwoon University, Seoul 139-701, Republic of Korea}
\email{tkkim@kw.ac.kr}

\author{DAE SAN KIM }
\address{Department of Mathematics, Sogang University, Seoul 121-742, Republic of Korea}
\email{dskim@sogang.ac.kr*}


\begin{abstract}

Dowling constructed Dowling lattice $Q_{n}(G)$, for any finite set with $n$ elements and any finite multiplicative group $G$ of order $m$, which is a finite geometric lattice. He also defined the Whitney numbers of the first and second kinds for any finite geometric lattice. These numbers for the Dowling lattice $Q_{n}(G)$ are the Whitney numbers of the first kind $V_{m}(n,k)$ and those of the second kind $W_{m}(n,k)$, which are given by Stirling number-like relations. In this paper, by `degenerating' such relations we introduce the degenerate Whitney numbers of the first kind and those of the second kind and investigate, among other things, generating functions, recurrence relations and various explicit expressions for them. As further generalizations of the degenerate Whitney numbers of both kinds, we also consider the degenerate $r$-Whitney numbers of both kinds.
\end{abstract}

\maketitle

\section{Introduction}
Let $(L,\le)$ be a finite lattice. This means that it is a finite poset such that every pair $x,y$ of elements in $L$ has a supremum $x \vee y$ and an infimum $x \wedge y$. It has a greatest element $\hat{1}$ and a least element $ \hat{0}$. $L$ is said to satisfy the chain condition if all maximal chains in any interval $[x, y]=\left\{z \in L|x \le z \le y \right\}$ have the same length. If $L$ satisfies the chain condition, then the rank $\rho(x)$ of an element $x \in L$ is the length of a maximal chain in $[\hat{0},x]$. The rank of a finite lattice $L$ is the maximum of lengths of all chains in $L$. An atom of a finite lattice $L$ is an element of rank 1, and $L$ is said to be atomic if every element of $L$ is the join of some atoms. A finite lattice $L$ is semimodular if $x \vee y$ covers both $x$ and $y$ whenever both $x$ and $y$ cover $x \wedge y$. Here, for $x,y \in L$, $y$ is said to cover $x$ iff $x < y$ and $[x,y]=\left\{x,y\right\}$. For more details on lattices, one refers to [24].\\
\indent A finite lattice $L$ is geometric if it is a finite semimodular lattice which is also atomic. Dowling [8] constructed an important finite geometric lattice $Q_{n}(G)$ out of a finite set of $n$ elements and a finite group $G$ of order $m$, called Dowling lattice of rank $n$ over a finite group of order $m$. Let $X$ be a finite set of $n $ elements and let $G$ be a finite multiplicative group of order $m$. Let $Q_{n}'(G)$ be the set of all partial $G$-partitions of $X$. Here a partial $G$-partition $\alpha=\left\{a_{j}: A_{j} \rightarrow G| j=1,2,\dots,r \right\}$ is a collection of functions where $A_{j}$'s are nonempty disjoint subsets of $X$ with $\cup_{j=1}^{r}A_{j}\subseteq X.$ Let  $\beta=\left\{b_{l}: B_{l} \rightarrow G| l=1,2,\dots,s \right\}$ be another partial $G$-partition of $X$. Then we write $ \alpha \le \beta$ if, for every $b_{l} \in \beta$, there exist a nonempty subset $\alpha_{l}$ of $\alpha$ and elements $g_{j} \in G$ such that $b_{l}=\sum_{a_{j} \in \alpha_{l}}g_{j}a_{j}$. This means that $B_{l}=\cup_{a_j \in \alpha_{l}}A_{j}$, and the restriction of $b_{l}$ to $A_{j}$ is given by $g_{j}a_{j}$, namely $b_{l}(x_{i})=g_{j}a_{j}(x_{i})$, for all $x_{i} \in A_{j}$. Then $\le$ defines a preorder on $Q_{n}'(G)$, so that the order relation, defined by $\alpha \sim \beta$ iff $\alpha \le \beta$ and $\beta \le \alpha$, gives an equivalence relation on $Q_{n}'(G)$. Now, the Dowling lattice $Q_{n}(G)$ of rank $n$ over a finite group of order $m$ is the quotient poset $Q_{n}'(G)/\sim$, with the order relation $[\alpha] \le [\beta]$ iff $ \alpha \le \beta$, where $[\alpha]$ denotes the equivalence class containing $\alpha$. Then it was shown that $Q_{n}(G)$ is a finite geometric lattice of rank $n$ (see [8,Theorem 3]) which satisfies the chain condition (see [8 ,Theorem 1]). We let the interested reader refer to [8] for the details on the construction of $Q_{n}(G)$ and its many fascinating properties. \\
\indent For a finite geometric lattice $L$ of rank $n$, Dowling [8] defined the Whitney numbers $V_{L}(n,k)$ of the first kind by
\begin{displaymath}
V(n,k)=V_{L}(n,k)=\sum_{\rho(x)=n-k}\mu(\hat{0},x),
\end{displaymath}
and the Whitney numbers $W_{L}(n,k)$ of the second kind by
\begin{displaymath}
W(n,k)=W_{L}(n,k)=\sum_{\rho(x)=n-k}1,
\end{displaymath}
where $\rho$ is the rank function and $\mu:L \times L \rightarrow \mathbb{Z}$ is the M\"{o}bius function given by
\begin{align*}
&\mu(x,y)=0, \quad \mathrm{for \,\, all}\,\, x \nleq y\,\,\mathrm{in}\,\,L, \\
&\mu(x,x)=1,\quad \mathrm{for \,\, all}\,\, x \in L, \\
&\mu(x,y)=-\sum_{x \le z < y}\mu(x,z),\quad \mathrm{for \,\, all}\,\, x < y\,\,\mathrm{in}\,\,L.
\end{align*}
If $L$ is the Dowling lattice $Q_{n}(G)$ of rank $n$ over a finite group $G$ of order $m$, then the Whitney numbers of the first kind $V_{Q_{n}(G)}(n,k)$ and the Whitney numbers of the second kind $W_{Q_{n}(G)}(n,k)$ are respectively denoted by $V_{m}(n,k)$ and $W_{m}(n,k)$. These notations are justified, since the Whitney numbers of both kinds depend only on the order $m$ of $G$.
In Corollary 6.1 of [8], it was shown that, for any fixed group $G$ of order $m$, the Whitney numbers $V_{m}(n,k)$ and $W_{m}(n,k)$ satisfy the following Stirling number-like relations:
\begin{align}
x^{n}=\sum_{k=0}^{n}W_{m}(n,k)m^{k}\bigg(\frac{x-1}{m}\bigg)_{k},\label{0}\\
m^{n}\bigg(\frac{x-1}{m}\bigg)_{n}=\sum_{k=0}^{n}V_{m}(n,k)x^{k},\label{1}
\end{align}
where $(x)_{0}=1$,\,\, $(x)_{n}=x(x-1)\cdots(x-n+1), \quad (n \ge 1)$.\par
Equivalently, the relations \eqref{0} and \eqref{1} are respectively given by
\begin{align}
&(mx+1)^{n}=\sum_{k=0}^{n}W_{m}(n,k)m^{k}(x)_{k},\label{2}\\
&m^{n}(x)_{n}=\sum_{k=0}^{n}V_{m}(n,k)(mx+1)^{k},\quad(\mathrm{see}\ [3,4,8]).\label{3}
\end{align}
Either from \eqref{0} and \eqref{1} or from \eqref{2} and \eqref{3}, it is evident that the Whitney numbers satisfy the orthogonality realtions. \par
\indent The aim of this paper is to introduce the degenerate Whitney numbers $V_{m,\lambda}(n,k)$ of the first kind and the degenerate Whitney numbers $W_{m,\lambda}(n,k)$ of the second kind as degenerate versions of the Whitney numbers of both kinds of Dowling lattice $Q_{n}(G)$ of rank $n$ over a finite group $G$ of order $m$ and investigate some combinatorial properties of those numbers, including generating functions, recurrence relations and various explicit expressions. As further generalizations of the degenerate Whitney numbers of both kinds, the degenerate $r$-Whitney numbers of the first kind $V_{m,\lambda}^{(r)}(n,k)$ and those of the second $W_{m,\lambda}^{(r)}(n,k)$ are studied and their connections with the degenerate $r$-Stirling numbers are observed. As immediate applications, it is expected that these numbers or their polynomial extensions will appear in the expressions of the probability distributions of appropriate random variables. Here we remark that the study of degenerate versions of some special numbers and polynomials, which was initiated by Carlitz in [6,7], regained the interests of some mathematicians and yielded many interesting results. For some of these, one refers to [1,11-19,22 and references therein]. Here we would like to mention only two things. In [12], starting from the question that what if we replace the usual exponential function appearing in the generating functions of Sheffer sequences by the degenerate exponential function in \eqref{11}, we were led to introduce the notions of degenerate Sheffer sequences and $\lambda$-Sheffer sequences. Also, it is worthwhile to note that studying degenerate versions can be done not only for polynomials but also for transcendental functions like gamma functions (see [15]). \par
\indent The outline of this paper is as follows. In Section 1, we recall some definitions, including Dowling polynomials, Tanny-Dowling polynomials, Stirling numbers of both kinds and Bell polynomials. And then we recall degenerate exponential functions, degenerate logarithms, degenerate Stirling numbers of both kinds, degenerate Bell polynomials and higher-order degenerate Bernoulli and Euler polynomials. Section 2 is the main results of this paper. In Section 2, we introduce the degenerate Whitney numbers of the second kind $W_{m,\lambda}(n,k)$ (see \eqref{21}) and its polynomial extensions, namely the degenerate Dowling polynomials $D_{m,\lambda}(n,x)$ (see \eqref{26}), and derive, among other things, their generating functions and explicit expressions for them when $m=1$, and their recurrence relations. In addition, we deduce various explicit expressions for the degenerate Whitney numbers of the second kind and a Dobinski-like formula for the degenerate Dowling polynomials. Then we introduce the degenerate Whitney numbers of the first kind $V_{m,\lambda}(n,k)$ by considering the inversion formula of \eqref{31}. We derive the generating function, a recurrence relation and various explicit expressions for them. Further, we introduce degenerate Tanny-Dowling polynomials and get the generating function of them. In Section 3, as further generalizations of the degenerate Whitney numbers of both kinds, we introduce the degenerate $r$-Whitney numbers $W_{m,\lambda}^{(r)}(n,k)$ of the second kind (see \eqref{65}) and the degenerate $r$-Whitney numbers $V_{m,\lambda}^{(r)}(n,k)$ of the first kind (see \eqref{64}), for any positive integer $r$ and derive their generating functions. Moreover, for $m=1$, we note the degenerate $r$-Whitney numbers of the second kind and those of the first kind are respectively given by the degenerate $r$-Stirling numbers of the second kind and those of the second kind. In the final Section 4, we get expressions for the higher-order degenerate Bernoulli and Euler numbers in terms of the degenerate Stirling numbers of both kinds.

For $n\ge 0$, Dowling polynomials are defined by
\begin{equation}
D_{m}(n,x)=\sum_{k=0}^{n}W_{m}(n,k)x^{k},\quad(\mathrm{see}\ [1,4]).\label{4}	
\end{equation}
When $x=1$, $D_{m}(n)=D_{m}(n,1)$ are called the Dowling numbers. \par
As is known, Tanny-Dowling polynomials are given by
\begin{align}
\mathcal{F}_{m}(n,x)\ &=\ \int_{0}^{+\infty}D_{m}(n,\mathcal{E}x)e^{-\mathcal{E}}d\mathcal{E}\label{5}	\\
&=\ \sum_{k=0}^{n}k!W_{m}(n,k)x^{k},\quad(\mathrm{see}\ [4]).\nonumber
\end{align}
The Stirling numbers of the first kind are defined by
\begin{equation}
(x)_{n}=\sum_{k=0}^{m}S_{1}(n,k)x^{k},\quad(n\ge 0),
\quad(\mathrm{see}\ [1,2,5,9,10,21]),\label{6}
\end{equation}
and the Stirling numbers of the second kind are given by
\begin{equation}
x^{k}=\sum_{k=0}^{n}S_{2}(n,k)(x)_{k},\quad(n\ge 0),\quad(\mathrm{see}\ [2,5,9,21-23]).\label{7}
\end{equation}
The Bell polynomials are defined by
\begin{equation}
\mathrm{Bel}_{n}(x)\ =\ \sum_{k=0}^{n}S_{2}(n,k)x^{k},\quad(n\ge 0),\quad(\mathrm{see}\ [16]).\label{8}
\end{equation}
The Fubini polynomials are given by
\begin{align}
\sum_{n=0}^{\infty}F_{n}(x)\frac{t^{n}}{n!}\ &=\ \frac{1}{1-x(e^{t}-1)}\ =\ \int_{0}^{+\infty}e^{\mathcal{E}x(e^{t}-1)}e^{-\mathcal{E}}d\mathcal{E}\label{9}\\
&=\ \sum_{n=0}^{\infty}\int_{0}^{+\infty}\mathrm{Bel}_{n}(\mathcal{E}x)e^{-\mathcal{E}}d\mathcal{E}\frac{t^{n}}{n!}.\nonumber	
\end{align}
Thus, by \eqref{8}, we get
\begin{equation}
F_{n}(x)=\sum_{k=0}^{n}k!S_{2}(n,k)x^{k},\quad(n\ge 0),\quad(\mathrm{see}\ [10,13]).\label{10}
\end{equation}
For any nonzero $\lambda\in\mathbb{R}$, the degenerate exponential functions are defined by
\begin{equation}
e_{\lambda}^{x}(t)\ =\ (1+\lambda t)^{\frac{x}{\lambda}}\ =\ \sum_{n=0}^{\infty}\frac{(x)_{n,\lambda}}{n!}{t^{n}},\quad(\mathrm{see}\ [11,13-15]),\label{11}
\end{equation}
where $(x)_{0,\lambda}=1,\ (x)_{n,\lambda}=x(x-\lambda)\cdots(x-(n-1)\lambda)$, $(n\ge 1)$. \par
When $x=1$, we denote $e_{\lambda}^{1}(t)$ simply by $e_{\lambda}(t)$.
Let $\log_{\lambda}(t)$ be the compositional inverse of $e_{\lambda}(t)$, called the degenerate logarithms. Then we have
\begin{equation}
\log_{\lambda}t=\frac{1}{\lambda}(t^{\lambda}-1),\quad \log_{\lambda}(1+t)=\sum_{n=1}^{\infty}
\lambda^{n-1}(1)_{n,1/\lambda}\frac{t^{n}}{n!},\quad(\mathrm{see}\ [11]).\label{12}
\end{equation}
In [11], the degenerate Stirling numbers of the first kind are defined by
\begin{equation}
(x)_{n}=\sum_{l=0}^{n}S_{1,\lambda}(n,l)(x)_{l,\lambda},\quad(n\ge 0). \label{13}	
\end{equation}
As an inversion formula of \eqref{13}, the degenerate Stirling numbers of the second kind are given by
\begin{equation}
(x)_{n,\lambda}=\sum_{k=0}^{n}S_{2,\lambda}(n,k)(x)_{k},\quad(n\ge 0),\quad(\mathrm{see}\ [11]).\label{14}
\end{equation}
The degenerate Bell polynomials are introduced by Kim-Kim-Dolgy as
\begin{equation}
e^{x(e_{\lambda}(t)-1)}=\sum_{n=0}^{\infty}\mathrm{Bel}_{n,\lambda}(x)\frac{t^{n}}{n!},\quad(\mathrm{see}\ [16]).\label{15}	
\end{equation}
Thus, we note that
\begin{equation}
\mathrm{Bel}_{n,\lambda}(x)=	\sum_{k=0}^{n}S_{2,\lambda}(n,k)x^{k},\quad(n\ge 0),\quad(\mathrm{see}\ [16]).\label{16}
\end{equation}
When $x=1$, we denote $\mathrm{Bel}_{n,\lambda}(1)$ simply by $\mathrm{Bel}_{n,\lambda}$. \par
From \eqref{13} and \eqref{14}, we note that
\begin{equation}
\frac{1}{k!}\big(e_{\lambda}(t)-1\big)^{k}=\sum_{n=k}^{\infty}S_{2,\lambda}(n,k)\frac{t^{n}}{n!},\label{17}
\end{equation}
and
\begin{equation}
\frac{1}{k!}\big(\log_{\lambda}(1+t)\big)^{k}=\sum_{n=k}^{\infty}S_{1,\lambda}(n,k)\frac{t^{n}}{n!},\quad(\mathrm{see}\ [11]).\label{18}	
\end{equation}
In [6,7], Carlitz considered the degenerate Bernoulli polynomials of order $r$ and the degenerate Euler polynomials of order $r$, respectively given by
\begin{equation}
\bigg(\frac{t}{e_{\lambda}(t)-1}\bigg)^{r}e_{\lambda}^{x}(t)	=\sum_{n=0}^{\infty}\beta_{n,\lambda}^{(r)}(x)\frac{t^{n}}{n!}, \label{19}
\end{equation}
and
\begin{equation}
\bigg(\frac{2}{e_{\lambda}(t)+1}\bigg)^{r}e_{\lambda}^{x}(t)	=\sum_{n=0}^{\infty}\mathcal{E}_{n,\lambda}^{(r)}(x)\frac{t^{n}}{n!}. \label{20}
\end{equation}
Note that $\displaystyle\lim_{\lambda\rightarrow 0}\beta_{n,\lambda}(x)=B_{n}^{(r)}(x) \displaystyle$ and $\displaystyle \lim_{\lambda\rightarrow 0}\mathcal{E}_{n,\lambda}^{(r)}(x)=E_{n}^{(r)}(x)$, where $B_{n}^{(r)}(x)$ are the Bernoulli polynomials of order $r$, and $E_{n}^{(r)}(x)$ are the Euler polynomials of order $r$.
When $r=1$, $\beta_{n,\lambda}^{(1)}(x)=\beta_{n,\lambda}(x)$, and $\mathcal{E}_{n,\lambda}^{(1)}(x)=\mathcal{E}_{n,\lambda}(x)$ are respectively called the degenerate Bernoulli polynomials and the degenerate Euler polynomials. For $x=0$, we let $\beta_{n,\lambda}^{(r)}(0)=\beta_{n,\lambda}^{(r)},\ \mathcal{E}_{n,\lambda}^{(r)}(0)=\mathcal{E}_{n,\lambda}^{(r)},\ \beta_{n,\lambda}(0)=\beta_{n,\lambda}$,\ and $\mathcal{E}_{n,\lambda}(0)=\mathcal{E}_{n,\lambda}$. \par

\section{degenerate Whitney numbers of Dowling lattices}
In light of \eqref{2}, we consider the degenerate Whitney numbers of the second kind given by
\begin{equation}
(mx+1)_{n,\lambda}=\sum_{k=0}^{n}W_{m,\lambda}(n,k)m^{k}(x)_{k},\quad(n\ge 0).\label{21}
\end{equation}
Now, we observe that
\begin{align}
e_{\lambda}^{mx+1}(t)\ &=\ \sum_{m=0}^{\infty}(mx+1)_{n,\lambda}\frac{t^{n}}{n!}\ =\ \sum_{n=0}^{\infty}\bigg(\sum_{k=0}^{n}W_{m,\lambda}(n,k)m^{k}(x)_{k}\bigg)\frac{t^{n}}{n!} \label{22}\\
&=\ \sum_{k=0}^{\infty}\bigg(\sum_{n=k}^{\infty}W_{m,\lambda}(n,k)\frac{t^{n}}{n!}\bigg)m^{k}(x)_{k}.\nonumber
\end{align}
On the other hand,
\begin{align}
e_{\lambda}^{mx+1}(t)\ &=\ e_{\lambda}(t)\cdot e_{\lambda}^{mx}(t)\ =\ e_{\lambda}(t)\cdot\big(e_{\lambda}^{m}(t)-1+1\big)^{x} \label{23} \\
&=\ e_{\lambda}(t)\sum_{k=0}^{\infty}\frac{1}{k!}\big(e_{\lambda}^{m}(t)-1\big)^{k}(x)_{k}.\nonumber
\end{align}
Therefore, by \eqref{22} and \eqref{23}, we obtain the following theorem.
\begin{theorem}
	For $k\ge 0$, we have
	\begin{displaymath}
		e_{\lambda}(t)\frac{1}{k!}\bigg(\frac{e_{\lambda}^{m}(t)-1}{m}\bigg)^{k}=\sum_{n=k}^{\infty}W_{m,\lambda}(n,k)\frac{t^{n}}{n!}.
	\end{displaymath}
\end{theorem}
In [11], we easily get
\begin{align*}
	&\sum_{n=k}^{\infty}S_{2,\lambda}(n+1,k+1)\frac{t^{n}}{n!}\ =\ \frac{d}{dt}\bigg(\sum_{n=k}^{\infty}S_{2,\lambda}(n+1,k+1)\frac{t^{n+1}}{(n+1)!}\bigg) \\
	&\ =\ \frac{d}{dt}\bigg(\frac{1}{(k+1)!}\big(e_{\lambda}(t)-1\big)^{k+1}\bigg)\ =\ \frac{1}{k!}\big(e_{\lambda}(t)-1\big)^{k}e_{\lambda}(t)\frac{1}{1+\lambda t}.
\end{align*}
Thus, we note that
\begin{align}
\frac{1}{k!}\big(e_{\lambda}(t)-1\big)^{k}e_{\lambda}(t)\ &=\ \sum_{n=k}^{\infty}S_{2,\lambda}(n+1,k+1)\frac{t^{n}}{n!}\big(1+\lambda	t\big) \label{24} \\
&=\ \sum_{n=k}^{\infty}\big(S_{2,\lambda}(n+1,k+1)+\lambda nS_{2,\lambda}(n,k+1)\big)\frac{t^{n}}{n!}. \nonumber
\end{align}
From Theorem 1 and \eqref{24}, we have
\begin{equation}
W_{1,\lambda}(n,k)=S_{2,\lambda}(n+1,k+1)+\lambda nS_{2,\lambda}(n,k+1), \label{25}
\end{equation}
where $n,k\ge 0$ with $n\ge k$.
\begin{corollary}
	For $n,k\ge 0$ with $n\ge k$, we have
\begin{displaymath}
	W_{1,\lambda}(n,k)=S_{2,\lambda}(n+1,k+1)+\lambda nS_{2,\lambda}(n,k+1).
\end{displaymath}	
\end{corollary}
In view of \eqref{4}, we define the degenerate Dowling polynomials as
\begin{equation}
D_{m,\lambda}(n,x)=\sum_{k=0}^{n}W_{m,\lambda}(n,k)x^{k},\quad(n\ge 0). \label{26}	
\end{equation}
When $x=1$, $D_{m,\lambda}(n,1)=D_{m,\lambda}(n)$ are called the degenerate Dowling numbers. \par
Then, by Theorem 1 and \eqref{26}, we get
\begin{align}
&\sum_{n=0}^{\infty}D_{m,\lambda}(n,x)\frac{t^{n}}{n!}\ =\ \sum_{n=0}^{\infty}\bigg(\sum_{k=0}^{n}W_{m,\lambda}(n,k)x^{k}\bigg)\frac{t^{n}}{n!}\label{27}\\
&=\ \sum_{k=0}^{\infty}x^{k}\sum_{n=k}^{\infty}W_{m,\lambda}(n,k)\frac{t^{n}}{n!}\ =\ e_{\lambda}(t)\sum_{k=0}^{\infty}x^{k}\frac{1}{k!}\bigg(\frac{e_{\lambda}^{m}(t)-1}{m}\bigg)^{k}\nonumber \\
&=	\ e_{\lambda}(t)e^{x(\frac{e_{\lambda}^{m}(t)-1}{m})}.\nonumber
\end{align}
Therefore, we obtain the generating function of degenerate Dowling polynomials.
\begin{theorem}
For $m\in\mathbb{N}$, we have
\begin{displaymath}
e_{\lambda}(t)e^{x(\frac{e_{\lambda}^{m}(t)-1}{m})}= \sum_{n=0}^{\infty}D_{m,\lambda}(n,x)\frac{t^{n}}{n!}.
\end{displaymath}
\end{theorem}
From \eqref{15}, we can easily derive the following equation
\begin{align}
&\sum_{n=0}^{\infty}\mathrm{Bel}_{n+1,\lambda}\frac{t^{n}}{n!}\ =\ \frac{d}{dt}e^{(e_{\lambda}(t)-1)}\ =\ \frac{d}{dt}\sum_{n=1}^{\infty}\bigg(\sum_{k=1}^{n}S_{2,\lambda}(n,k)\bigg)\frac{t^{n}}{n!} \label{28} \\
&=\ \sum_{n=0}^{\infty}\bigg(\sum_{k=1}^{n+1}S_{2,\lambda}(n+1,k)\bigg)\frac{t^{n}}{n!}\ =\ \sum_{n=0}^{\infty}\bigg(\sum_{k=0}^{n}S_{2,\lambda}(n+1,k+1)\bigg)\frac{t^{n}}{n!}. \nonumber	
\end{align}
Thus, we have
\begin{align}
&\mathrm{Bel}_{n+1,\lambda}=\sum_{k=0}^{n}S_{2,\lambda}(n+1,k+1),\label{29}	\\
&\mathrm{Bel}_{n+1,\lambda}(x)=x\sum_{k=0}^{n}S_{2,\lambda}(n+1,k+1)x^{k},\quad(n\ge 0).\nonumber
\end{align}
Here we note that \eqref{29} corresponds to the fact that $S_{2,\lambda}(n+1,0)=0$, for $n \ge 0$.

Using \eqref{25}, we have
\begin{align}
D_{1,\lambda}(n)\ &=\ \sum_{k=0}^{n}W_{1,\lambda}(n,k)\ =\ \sum_{k=0}^{n}\big\{S_{2,\lambda}(n+1,k+1)+n\lambda S_{2,\lambda}(n,k+1)\big\}\label{30} \\
&=\ \mathrm{Bel}_{n+1,\lambda}+n\lambda \sum_{k=1}^{n}S_{2,\lambda}(n,k)\ =\ \mathrm{Bel}_{n+1,\lambda}+n \lambda\sum_{k=0}^{n}S_{2,\lambda}(n,k) \nonumber \\
&=\ \mathrm{Bel}_{n+1,\lambda}+n\lambda\mathrm{Bel}_{n,\lambda}. \nonumber	
\end{align}
\begin{corollary}
	For $n\ge 0$, we have
	\begin{displaymath}
		D_{1,\lambda}(n)= \mathrm{Bel}_{n+1,\lambda}+n\lambda\mathrm{Bel}_{n,\lambda}.
	\end{displaymath}
\end{corollary}
As an inversion formula of \eqref{21}, we define the degenerate Whitney numbers of the first kind as
\begin{equation}
m^{n}(x)_{n}=\sum_{k=0}^{n}V_{m,\lambda}(n,k)(mx+1)_{k,\lambda},\quad(n\ge 0). \label{31}
\end{equation}
From \eqref{31}, we note that
\begin{equation}
	\begin{aligned}
	(1+mt)^{x}	\ &=\ \sum_{n=0}^{\infty}(x)_{n}m^{n}\frac{t^{n}}{n!}\ =\ \sum_{n=0}^{\infty}\bigg(\sum_{k=0}^{n}V_{m,\lambda}(n,k)(mx+1)_{k,\lambda}\bigg)\frac{t^{n}}{n!} \\
&=\ \sum_{k=0}^{\infty}\bigg(\sum_{n=k}^{\infty}V_{m,\lambda}(n,k)\frac{t^{n}}{n!}\bigg)(mx+1)_{k,\lambda}.	\end{aligned}\label{32}
\end{equation}
On the other hand, with the help of \eqref{12} we have
\begin{align}
&(1+mt)^{x}=(1+mt)^{x+\frac{1}{m}}(1+mt)^{-\frac{1}{m}}=(1+mt)^{\frac{mx+1}{m}} (1+mt)^{-\frac{1}{m}}\label{33}\\
&\ =\ e_{\lambda}^{mx+1}\big(\log_{\lambda}(1+mt)^{\frac{1}{m}}\big)(1+mt)^{-\frac{1}{m}} \nonumber \\
&\ =\ \sum_{k=0}^{\infty}\frac{1}{k!}\bigg(\log_{\lambda}(1+mt)^{\frac{1}{m}}\bigg)^{k}(1+mt)^{-\frac{1}{m}}(mx+1)_{k,\lambda}.\nonumber	
\end{align}
From \eqref{11}, we note that
\begin{displaymath}
	e_{m}(t)=(1+mt)^{\frac{1}{m}},\quad e_{m}^{-1}(t)=(1+mt)^{-\frac{1}{m}}.
\end{displaymath}
Therefore, by \eqref{32} and \eqref{33}, we obtain the generating function of the degenerate Whitney numbers of the first kind.
\begin{theorem}
	For $k\ge 0$, we have
	\begin{displaymath}
		\frac{1}{k!}\big(\log_{\lambda} e_{m}(t)\big)^{k}e_{m}^{-1}(t)=\sum_{n=k}^{\infty}V_{m,\lambda}(n,k)\frac{t^{n}}{n!}.
	\end{displaymath}
\end{theorem}
By \eqref{21}, we get
\begin{align}
&\quad \sum_{k=0}^{n}W_{m,\lambda}(n,k)m^{k}(x)_{k}\ =\ (mx+1)_{n,\lambda}\ =\ (mx+1)_{n-1,\lambda}\big(mx+1-\lambda(n-1)\big)\label{34}	\\
&=\ \sum_{k=0}^{n-1}W_{m,\lambda}(n-1,k)m^{k}(x)_{k}\big(m(x-k)+mk+1-\lambda(n-1)\big) \nonumber \\
&=\ \sum_{k=0}^{n-1}W_{m,\lambda}(n-1,k)m^{k+1}(x)_{k+1}+\sum_{k=0}^{n-1}W_{m,\lambda}(n-1,k)\big(mk+1-\lambda(n-1)\big)m^{k}(x)_{k}\nonumber\\
&=\ \sum_{k=0}^{n}\big\{W_{m,\lambda}(n-1,k-1)+(mk+1)W_{m,\lambda}(n-1,k)-\lambda(n-1)W_{m,\lambda}(n-1,k)\big\}m^{k}(x)_{k}.\nonumber
\end{align}
Therefore, by \eqref{34}, we obtain the following theorem.
\begin{theorem}
	For $n,k\ge 0$ with $n\ge k$, we have
	\begin{displaymath}
		W_{m,\lambda}(n,k)= W_{m,\lambda}(n-1,k-1)+(mk+1)W_{m,\lambda}(n-1,k)-\lambda(n-1)W_{m,\lambda}(n-1,k).
	\end{displaymath}
\end{theorem}
From \eqref{31}, we have
\begin{align}
	&\sum_{k=0}^{n}V_{m,\lambda}(n,k)(mx+1)_{k,\lambda}\ =\ m^{n}(x)_{n}\ =\ m(x-n+1)m^{n-1}(x)_{n-1} \label{35} \\
	&=\ \sum_{k=0}^{n-1}V_{m,\lambda}(n-1,k)(mx+1)_{k,\lambda}(mx+1-k\lambda+k\lambda+m-nm-1)\nonumber \\
	&=\ \sum_{k=0}^{n-1}V_{m,\lambda}(n-1,k)(mx+1)_{k+1,\lambda}+\sum_{k=0}^{n-1}(k\lambda+m-nm-1)V_{m,\lambda}(n-1,k)(mx+1)_{k,\lambda}\nonumber \\
	&=\ \sum_{k=0}^{n}\big\{V_{m,\lambda}(n-1,k-1)+(k\lambda+m-nm-1)V_{m,\lambda}(n-1,k)\big\}(mx+1)_{k,\lambda}.\nonumber
\end{align}
	Therefore, by comparing the coefficients on both sides of \eqref{35}, we obtain the following theorem.
\begin{theorem}
For $n,k\ge 0$ with $n\ge k$, we have
\begin{displaymath}
V_{m,\lambda}(n,k)=(m-nm-1)V_{m,\lambda}(n-1,k)+ V_{m,\lambda}(n-1,k-1)+k\lambda V_{m,\lambda}(n-1,k).
\end{displaymath}
\end{theorem}
\noindent \emph{Remark.} From Theorem 7 and \eqref{31}, we obtain
\begin{equation}
V_{m,\lambda}(n,n)=V_{m,\lambda}(n-1,n-1)=\cdots=V_{m,\lambda}(0,0)=1.\label{35-1}
\end{equation}
By Theorem 5, we get
\begin{align}
&\sum_{n=k}^{\infty}V_{m,\lambda}(n,k)\frac{t^{n}}{n!}\ =\ \frac{1}{k!}\big(\log_{\lambda}(1+mt)^{\frac{1}{m}}\big)^{k}(1+mt)^{-\frac{1}{m}} \label{36} \\
&=\ \frac{1}{k!}\big(\log_{\lambda}(1+mt)^{\frac{1}{m}}\big)^{k}e_{\lambda}^{-1}\big(\log_{\lambda}(1+mt)^{\frac{1}{m}}\big)\nonumber \\
&=\ \sum_{l=0}^{\infty}\frac{(-1)_{l,\lambda}}{l!}\big(\log_{\lambda}(1+mt)^{\frac{1}{m}}\big)^{l}\frac{1}{k!}\big(\log_{\lambda}(1+mt)^{\frac{1}{m}}\big)^{k}\nonumber\\
&=\ \sum_{l=0}^{\infty}\frac{(-1)^{l}\langle l\rangle_{l,\lambda}}{l!k!}\frac{(k+l)!}{(k+l)!}\big(\log_{\lambda}(1+mt)^{\frac{1}{m}}\big)^{k+l}\nonumber \\
&=\ \sum_{l=0}^{\infty}(-1)^{l}\binom{k+l}{k}\langle 1\rangle_{l,\lambda}\frac{1}{(k+l)!}\Big(\log_{\lambda}\big((1+mt)^{\frac{1}{m}}-1+1)\Big)^{l+k}\nonumber\\
&=\ \sum_{l=0}^{\infty}\binom{k+l}{k}(-1)^{l}\langle 1\rangle_{l,\lambda}\sum_{i=k+l}^{\infty}S_{1,\lambda}(i,k+l)\frac{1}{i!}\big((1+mt)^{\frac{1}{m}}-1\big)^{i}\nonumber\\
&=\ \sum_{l=k}^{\infty}\binom{l}{k}(-1)^{l-k}\langle 1\rangle_{l-k,\lambda}\sum_{i=l}^{\infty}S_{1,\lambda}(i,l)\frac{1}{i!}\Big(e^{\frac{1}{m}\log(1+mt)}-1\Big)^{i}\nonumber\\
&=\ \sum_{i=k}^{\infty}\bigg(\sum_{l=k}^{i}\binom{l}{k}(-1)^{l-k}\langle 1\rangle_{l-k,\lambda}S_{1,\lambda}(i,l)\bigg)\sum_{j=i}^{\infty}S_{2}(j,i)m^{-j}\frac{1}{j!}\big(\log(1+mt)\big)^{j} \nonumber \\
&=\ \sum_{j=k}^{\infty}\bigg(\sum_{i=k}^{j}\sum_{l=k}^{i}\binom{l}{k}(-1)^{l-k}\langle 1\rangle_{l-k,\lambda}S_{1,\lambda}(i,l)S_{2}(j,i)m^{-j}\bigg)\sum_{n=j}^{\infty}S_{1}(n,j)\frac{m^{n}t^{n}}{n!}\nonumber \\
&=\ \sum_{n=k}^{\infty}\bigg(\sum_{j=k}^{n}\sum_{i=k}^{j}\sum_{l=k}^{i}\binom{l}{k}(-1)^{l-k}\langle 1\rangle_{l-k,\lambda}S_{1,\lambda}(i,l)S_{2}(j,i)S_{1}(n,j)m^{n-j}\bigg)\frac{t^{n}}{n!}, \nonumber
\end{align}
where $\langle x\rangle_{0,\lambda}=1$, $\langle x\rangle_{n,\lambda}=x(x+\lambda)\cdots(x+(n-1)\lambda),\ (n\ge 1)$. \par
By comparing the coefficients on both sides of \eqref{36}, we have
\begin{align}
V_{m,\lambda}(n,k)\ &=\ \sum_{j=k}^{n}\sum_{i=k}^{j}\sum_{l=k}^{i}\binom{l}{k}(-1)^{l-k}\langle 1\rangle_{l-k,\lambda}S_{1,\lambda}(i,l)S_{2}(j,i)S_{1}(n,j)m^{n-j}\label{37}\\
&=\ \sum_{j=k}^{n}\sum_{l=k}^{j}\sum_{i=l}^{j}\binom{l}{k}(-1)^{l-k}\langle 1\rangle_{l-k,\lambda}S_{1,\lambda}(i,l)S_{2}(j,i)S_{1}(n,j)m^{n-j}.\nonumber
\end{align}
Therefore, we obtain the following theorem.
\begin{theorem}
	For $n,k\ge 0$ with $n\ge k$, we have
	\begin{displaymath}
		V_{m,\lambda}(n,k)=\sum_{j=k}^{n}\sum_{l=k}^{j}\binom{l}{k}(-1)^{l-k}\langle 1\rangle_{l-k,\lambda}\bigg(\sum_{i=l}^{j}S_{1,\lambda}(i,l)S_{2}(j,i)\bigg)S_{1}(n,j)m^{n-j}.
	\end{displaymath}
\end{theorem}
We note that, from (3), (31) and Theorem 8, we have
\begin{align*}
	\lim_{\lambda\rightarrow 0}V_{m,\lambda}(n,k)\ &=\ V_{m}(n,k) \\
	&=\ \sum_{l=k}^{n}\binom{l}{k}(-1)^{l-k}S_{1}(n,l)m^{n-l}.
\end{align*}
In view of \eqref{5}, we may consider the degenerate Tanny-Dowling polynomials given by
\begin{align}
\mathcal{F}_{m,\lambda}(n,x)\ &=\ \int_{0}^{+\infty}D_{m,\lambda}(n,x\mathcal{E})e^{-\mathcal{E}}d\mathcal{E} \label{38} \\
&=\ \sum_{k=0}^{n}W_{m,\lambda}(n,k)x^{k}\int_{0}^{+\infty}e^{-\mathcal{E}}\mathcal{E}^{k}d\mathcal{E} \nonumber \\
&=\sum_{k=0}^{n}W_{m,\lambda}(n,k)k!x^k. \nonumber
\end{align}
From Theorem 3, we have
\begin{align*}
	\sum_{n=0}^{\infty}\mathcal{F}_{m,\lambda}(n,x)\frac{t^{n}}{n!}\ &=\ \sum_{n=0}^{\infty}\int_{0}^{+\infty}D_{m,\lambda}(n,x\mathcal{E})e^{-\mathcal{E}}d\mathcal{E}\frac{t^{n}}{n!} \\
	&=\ \int_{0}^{+\infty}e^{-\mathcal{E}}\sum_{n=0}^{\infty}D_{m,\lambda}(n,x\mathcal{E})\frac{t^{n}}{n!}d\mathcal{E}\\
	&=\ e_{\lambda}(t)\int_{0}^{\infty}e^{-\mathcal{E}(1-x(\frac{e_{\lambda}^{m}(t)-1}{m}))}d\mathcal{E}\\
 &=\ e_{\lambda}(t)\frac{1}{1-x\big(\frac{e_{\lambda}^{m}(t)-1}{m}\big)}.
\end{align*}
Therefore, we obtain the following theorem.
\begin{theorem}
For $m\in\mathbb{N}$, we have
\begin{displaymath}
e_{\lambda}(t)\frac{1}{1-x\big(\frac{e_{\lambda}^{m}(t)-1}{m}\big)}=	\sum_{n=0}^{\infty}\mathcal{F}_{m,\lambda}(n,x)\frac{t^{n}}{n!}.
\end{displaymath}
\end{theorem}
From Theorem 3, we note that
\begin{align}
\sum_{n=0}^{\infty}D_{m,\lambda}(n,x)\frac{t^{n}}{n!}\ &=\ e_{\lambda}(t)e^{x(\frac{e_{\lambda}^{m}(t)-1}{m})}=e^{-\frac{x}{m}}\sum_{k=0}^{\infty}\frac{x^{k}e_{\lambda}^{mk+1}(t)}{m^{k}k!}\label{39} \\
&= e^{-\frac{x}{m}}\sum_{k=0}^{\infty}\frac{x^{k}}{k!m^{k}}\sum_{n=0}^{\infty}(mk+1)_{n,\lambda}\frac{t^{n}}{n!}\nonumber\\
&=\ \sum_{n=0}^{\infty}\bigg\{e^{-\frac{x}{m}}\sum_{k=0}^{\infty}\frac{x^{k}}{k!m^{k}}(mk+1)_{n,\lambda}\bigg\}\frac{t^{n}}{n!}.\nonumber	
\end{align}
Therefore, we obtain Dobinski-like formula for the degenerate Dowling polynomials.
\begin{theorem}
	For $n\ge 0$, we have
	\begin{displaymath}
		D_{m,\lambda}(n,x)=e^{-\frac{x}{m}}\sum_{k=0}^{\infty}\frac{x^{k}}{m^{k}k!}(mk+1)_{n,\lambda}.
	\end{displaymath}
	In particular,
	\begin{displaymath}
		D_{m,\lambda}(n)=e^{-\frac{1}{m}}\sum_{k=0}^{\infty}\frac{1}{m^{k}k!}(mk+1)_{n,\lambda}.
	\end{displaymath}	
\end{theorem}
For the next result, we observe from \eqref{15} that $\mathrm{Bel}_{n,\lambda}(x)=\frac{1}{e^x}\sum_{k=0}^{\infty}\frac{(k)_{n,\lambda}}{k!}x^k$. From Theorem 10 with $m=1$ and for $n \ge 0$, we have
\begin{align*}
D_{1,\lambda}(n,x)\ &=\ e^{-x}\sum_{k=0}^{\infty}\frac{x^{k}}{k!}(k+1)_{n,\lambda}\ =\ e^{-x}\sum_{k=1}^{\infty}\frac{x^{k-1}}{(k-1)!}(k)_{n,\lambda} \\
&=\ \frac{1}{xe^{x}}\sum_{k=1}^{\infty}\frac{(k)_{n,\lambda}}{k!}x^{k}(k-n\lambda+n\lambda) \\
&=\ \frac{1}{xe^{x}}\bigg(\sum_{k=1}^{\infty}\frac{(k)_{n+1,\lambda}}{k!}x^{k}+n\lambda\sum_{k=1}^{\infty}\frac{(k)_{n,\lambda}}{k!}x^{k}\bigg) \\
&=\ \frac{1}{x}\bigg(\frac{1}{e^{x}}\sum_{k=1}^{\infty}\frac{(k)_{n+1,\lambda}}{k!}x^{k}+n\lambda\frac{1}{e^{x}}\sum_{k=1}^{\infty}\frac{(k)_{n,\lambda}}{k!}x^{k}\bigg) \\
&=\ \frac{1}{x}\bigg(\frac{1}{e^{x}}\sum_{k=0}^{\infty}\frac{(k)_{n+1,\lambda}}{k!}x^{k}+n\lambda\frac{1}{e^{x}}\sum_{k=0}^{\infty}\frac{(k)_{n,\lambda}}{k!}x^{k}\bigg) \\
&=\frac{1}{x}\big(\mathrm{Bel}_{n+1,\lambda}(x)+n \lambda\mathrm{Bel}_{n,\lambda}(x)\big).
\end{align*}
Thus, we get the following corollary.
\begin{corollary}
	For $n\ge 0$, we have
	\begin{displaymath}
	xD_{1,\lambda}(n,x)=\mathrm{Bel}_{n+1,\lambda}(x)+n\lambda\mathrm{Bel}_{n,\lambda}(x).
\end{displaymath}
\end{corollary}
By Theorem 1, we get
\begin{align}
&\sum_{n=k}^{\infty}W_{m,\lambda}(n,k)\frac{t^{n}}{n!}\ =\ e_{\lambda}(t)\frac{1}{k!}\bigg(\frac{e_{\lambda}^{m}(t)-1}{m}\bigg)^{k} \label{40} \\
&=\ \frac{1}{k!m^{k}}\sum_{l=0}^{k}\binom{k}{l}(-1)^{k-l}e_{\lambda}^{lm+1}(t)=\frac{1}{k!m^{k}}\sum_{l=0}^{k}\binom{k}{l}(-1)^{k-l}\sum_{n=0}^{\infty}(lm+1)_{n,\lambda}\frac{t^{n}}{n!} \nonumber \\
&=\ \sum_{n=0}^{\infty}\bigg(\frac{1}{k!m^{k}}\sum_{l=0}^{k}\binom{k}{l}(-1)^{k-l}(lm+1)_{n,\lambda}\bigg)\frac{t^{n}}{n!}. \nonumber	
\end{align}
Comparing the coefficients on both sides of \eqref{40}, we obtain the following theorem.
\begin{theorem}
For $n,k\ge 0$, we have
\begin{displaymath}
\frac{1}{k!m^{k}}\sum_{l=0}^{k}\binom{k}{l}(-1)^{k-l}(lm+1)_{n,\lambda}=\left\{\begin{array}{cc}
W_{m,\lambda}(n,k), & \textrm{if $n\ge k,$}\\
0, & \textrm{otherwise.}
\end{array}\right.
\end{displaymath}
\end{theorem}
From Theorem 12 with $m=1$ and Corollary 2,  for $n\ge k\ge 0$, we have
\begin{align*}
\frac{1}{k!}\sum_{l=0}^{k}\binom{k}{l}(-1)^{k-l}(l+1)_{n,\lambda}\ &=\ \frac{1}{k!}\sum_{l=0}^{k}\binom{k}{l}(-1)^{l}(k-l+1)_{n,\lambda}\ =\ W_{1,\lambda}(n,k) \\
&=\ S_{2,\lambda}(n+1,k+1)+n\lambda S_{2,\lambda}(n,k+1).
\end{align*}
Let $n,k\ge 0$ with $n\ge k$. Then, by Theorem 12, we get
\begin{align}
W_{m,\lambda}(n,k)\ &=\ \frac{1}{m^{k}k!}\sum_{l=0}^{k}\binom{k}{l}(-1)^{l}\big(m(k-l)+1\big)_{n,\lambda} \label{41} \\
&=\ \frac{1}{m^{k}k!}\sum_{l=0}^{k}(-1)^{l}\binom{k}{l}\sum_{i=0}^{n}\binom{n}{i}\big(m(k-l)\big)_{i,\lambda}(1)_{n-i,\lambda}\nonumber 	\\
&=\ \frac{1}{m^{k}k!}\sum_{l=0}^{k}(-1)^{l}\binom{k}{l}\sum_{i=0}^{n}\binom{n}{i}m^{i}(k-l)_{i,\frac{\lambda}{m}}(1)_{n-i,\lambda}\nonumber \\
&=\ \sum_{i=0}^{n}\binom{n}{i}m^{i-k}(1)_{n-i,\lambda}\frac{1}{k!}\sum_{l=0}^{k}(-1)^{l}\binom{k}{l}(k-l)_{i,\frac{\lambda}{m}} \nonumber \\
&=\ \sum_{i=k}^{n}\binom{n}{i}m^{i-k}(1)_{n-i,\lambda}\frac{1}{k!}\sum_{l=0}^{k}(-1)^{l}\binom{k}{l}(k-l)_{i,\frac{\lambda}{m}} \nonumber \\
&=\ \sum_{i=k}^{n}\binom{n}{i}m^{i-k}(1)_{n-i,\lambda}S_{2,\frac{\lambda}{m}}(i,k).\nonumber
\end{align}
Therefore, by \eqref{41}, we obtain the following theorem.
\begin{theorem}
	For $n,k\ge 0$ with $n\ge k$, we have
	\begin{displaymath}
		W_{m,\lambda}(n,k)= \sum_{i=k}^{n}\binom{n}{i}m^{i-k}(1)_{n-i,\lambda}S_{2,\frac{\lambda}{m}}(i,k).
	\end{displaymath}
\end{theorem}
Note that
\begin{displaymath}
	W_{1,\lambda}(n,k)=S_{2,\lambda}(n+1,k+1)+n\lambda
 S_{2,\lambda}(n,k+1)=\sum_{i=k}^{n}\binom{n}{i}(1)_{n-i,\lambda}S_{2,\lambda}(i,k).
 \end{displaymath}
 Let $\triangle$ be the difference operator with $\triangle f(x)=f(x+1)-f(x)$. Then we easily get
 \begin{displaymath}
 	\triangle^{k}f(x)=\sum_{l=0}^{k}\binom{k}{l}(-1)^{k-l}f(x+l), (k\ge 0).
 \end{displaymath}
 Let us take $f(x)=(mx+1)_{n,\lambda}$. Then we have
\begin{equation}
 \sum_{l=0}^{k}\binom{k}{l}(-1)^{k-l}(lm+1)_{n,\lambda}=\triangle^{k}(mx+1)_{n,\lambda}\big|_{x=0}.\label{42}
\end{equation}
From Theorem 12 and \eqref{42}, for $n \ge k \ge 0$, we have
\begin{align}
k!m^{k}W_{m,\lambda}(n,k)\ &=\ \sum_{l=0}^{k}\binom{k}{l}(-1)^{k-l}(lm+1)_{n,\lambda} \nonumber \\
&=\ \triangle^{k}(mx+1)_{n,\lambda}\big|_{x=0}.\label{43}
\end{align}
Therefore, by \eqref{43}, we obtain the following theorem.
\begin{theorem}
For $n,k\ge 0$ with $n\ge k$, we have
\begin{displaymath}
W_{m,\lambda}(n,k)=\frac{1}{k!m^{k}}\triangle^{k}(mx+1)_{n,\lambda}\bigg|_{x=0}.
\end{displaymath}
\end{theorem}
For $k\ge 0$, we observe that
\begin{align}
&\frac{1}{k!}e_{\lambda}^{z}(t)\big(1-e_{\lambda}^{-1}(t)\big)^{k}\ =\ \frac{1}{k!}e_{\lambda}^{z-k}(t)\big(e_{\lambda}(t)-1\big)^{k} \label{44} \\
&\ =\ \sum_{l=0}^{\infty}(z-k)_{l,\lambda}\frac{t^{l}}{l!}\sum_{i=k}^{\infty}S_{2,\lambda}(i,k)\frac{t^{i}}{i!} \nonumber \\
&\ =\ \sum_{n=k}^{\infty}\bigg(\sum_{i=k}^{n}\binom{n}{i}S_{2,\lambda}(i,k)(z-k)_{n-i,\lambda}\bigg)\frac{t^{n}}{n!}\nonumber.
\end{align}
On the other hand,
\begin{align}
\frac{1}{k!}e_{\lambda}^{z}(t)\big(1-e_{\lambda}^{-1}(t)\big)^{k}\ &=\ \frac{1}{k!}\sum_{l=0}^{k}\binom{k}{l}(-1)^{l}e_{\lambda}^{z-l}(t) \label{45}\\
&=\ \sum_{n=0}^{\infty}\bigg(\frac{1}{k!}\sum_{l=0}^{k}\binom{k}{l}(-1)^{l}(z-l)_{n,\lambda}\bigg)\frac{t^{n}}{n!}.\nonumber	
\end{align}
Thus, by \eqref{44} and \eqref{45}, we get
\begin{equation}
\frac{1}{k!}\sum_{l=0}^{k}\binom{k}{l}(-1)^{l}(z-l)_{n,\lambda}=\left\{\begin{array}{ccc}
	\displaystyle\sum_{i=k}^{n}\binom{n}{i}S_{2,\lambda}(i,k)(z-k)_{n-i,\lambda}\displaystyle, & \textrm{if $n\ge k,$} \\
	0, & \textrm{otherwise.}
\end{array}\right.\label{46}
\end{equation}
Let us take $n=k$. Then, by \eqref{46}, we easily get
\begin{equation}
\sum_{l=0}^{k}\binom{k}{l}(-1)^{l}(z-l)_{k,\lambda}=k!,\label{47}
\end{equation}
where $z$ is any complex number. \par
Now, we are going to show \eqref{47} in another way.
For $n\ge 0$, by Theorem 12, we get
\begin{align}
W_{m,\lambda}(n,n)\ &=\ \frac{1}{m^{n}n!}\sum_{j=0}^{n}(-1)^{j}\binom{n}{j}m^{n}\bigg(n-j+\frac{1}{m}\bigg)_{n,\frac{\lambda}{m}}\label{48}	\\
&=\ \frac{1}{n!}\sum_{j=0}^{n}(-1)^{j}\binom{n}{j}\bigg(n-j+\frac{1}{m}\bigg)_{n,\frac{\lambda}{m}}.\nonumber
\end{align}
From \eqref{48} and \eqref{35-1}, we note that
\begin{equation}
1=\frac{1}{n!}\sum_{j=0}^{n}(-1)^{j}\binom{n}{j}\bigg(n-j+\frac{1}{m}\bigg)_{n,\frac{\lambda}{m}}.\label{49}	
\end{equation}
We observe that the right-hand side of \eqref{49} is independent of $m,\lambda$. So, replacing
$m$ by $\frac{1}{z-n}$ and $\lambda$ by $\frac{\lambda}{z-n}$, we obtain
\begin{equation}
1=\frac{1}{n!}\sum_{j=0}^{n}(-1)^{j}\binom{n}{j}(z-j)_{n,\lambda}.\label{50}	
\end{equation}
Therefore, by either \eqref{47} or \eqref{50}, we obtain the following lemma.
\begin{lemma}
	For $n\ge 0$ and any complex  number $z$, we have
	\begin{displaymath}
		\sum_{j=0}^{n}(-1)^{j}\binom{n}{j}(z-j)_{n,\lambda}=n!.
	\end{displaymath}
\end{lemma}
From Theorem 1, we note that
\begin{align}
&\sum_{n=k-1}^{\infty}W_{m,\lambda}(n+1,k)\frac{t^{n}}{n!}\ =\ \frac{d}{dt}\frac{1}{k!}e_{\lambda}(t)\bigg(\frac{e_{\lambda}^{m}(t)-1}{m}\bigg)^{k}\label{51} \\
&=\ \frac{e_{\lambda}^{1-\lambda}(t)}{k!} \bigg(\frac{e_{\lambda}^{m}(t)-1}{m}\bigg)^{k}+\frac{ke_{\lambda}(t)}{k!} \bigg(\frac{e_{\lambda}^{m}(t)-1}{m}\bigg)^{k-1}e_{\lambda}^{m-\lambda}(t)\nonumber \\
&=\ \bigg\{\frac{1}{k!}e_{\lambda}(t) \bigg(\frac{e_{\lambda}^{m}(t)-1}{m}\bigg)^{k}+\frac{1}{(k-1)!}e_{\lambda}(t) \bigg(\frac{e_{\lambda}^{m}(t)-1}{m}\bigg)^{k-1}e_{\lambda}^{m}(t)\bigg\}\frac{1}{1+\lambda t}\nonumber \\
&=\ \sum_{l=k-1}^{\infty}\bigg\{W_{m,\lambda}(l,k)+\sum_{i=k-1}^{l}W_{m,\lambda}(i,k-1)(m)_{l-i,\lambda}\binom{l}{i}\bigg\}\frac{t^{l}}{l!}\sum_{j=0}^{\infty}(-1)^{j}\lambda^{j}t^{j} \nonumber \\
&=\ \sum_{n=k-1}^{\infty}\bigg(n!\sum_{l=k-1}^{n}\bigg(W_{m,\lambda}(l,k)+\sum_{i=k-1}^{l}W_{m,\lambda}(i,k-1)\binom{l}{i}\bigg)\frac{(-\lambda)^{n-l}}{l!}\bigg)\frac{t^{n}}{n!}. \nonumber	
\end{align}
Therefore, we \eqref{51}, we obtain the following theorem.
\begin{theorem}
For $0\le k\le n$, we have the recursion formula :
\begin{displaymath}
W_{m,\lambda}(n+1,k)= n!\sum_{l=k-1}^{n}\bigg(W_{m,\lambda}(l,k)+\sum_{i=k-1}^{l}W_{m,\lambda}(i,k-1)\binom{l}{i}\bigg)\frac{(-\lambda)^{n-l}}{l!}.
\end{displaymath}
\end{theorem}
Note that
\begin{displaymath}
	W_{m}(n+1,k)=\lim_{\lambda\rightarrow 0}W_{m,\lambda}(n+1,k)=W_{m}(n,k)+\sum_{i=k-1}^{n}\binom{n}{i}W_{m}(i,k-1).
\end{displaymath}

From Theorem 3, we note that
\begin{align}
&\sum_{n=0}^{\infty}D_{m,\lambda}(n+1)\frac{t^{n}}{n!}\ =\ \frac{d}{dt}e_{\lambda}(t)e^{(\frac{e_{\lambda}^{m}(t)-1}{m})} \label{52} \\
&\ =\ e_{\lambda}^{1-\lambda}(t)e^{(\frac{e_{\lambda}^{m}(t)-1}{m})}+e_{\lambda}(t)e^{\frac{e_{\lambda}^{m}(t)-1}{m}}\cdot e_{\lambda}^{m-\lambda}(t)\nonumber \\
&=\ e_{\lambda}^{-\lambda}(t)\sum_{l=0}^{\infty}D_{m,\lambda}(l)\frac{t^{l}}{l!}+e_{\lambda}^{-\lambda}(t)\sum_{j=0}^{\infty}D_{m,\lambda}(j)\frac{t^{j}}{j!}\sum_{i=0}^{\infty}(m)_{i,\lambda}\frac{t^{i}}{i!}\nonumber
\end{align}
\begin{align*}
&=\ e_{\lambda}^{-\lambda}(t)\sum_{l=0}^{\infty}\bigg(D_{m,\lambda}(l)+\sum_{i=0}^{l}\binom{l}{i}D_{m,\lambda}(l-i)(m)_{i,\lambda}\bigg)\frac{t^{l}}{l!} \nonumber \\
&=\ \sum_{k=0}^{\infty}(-\lambda)^{k}t^{k}\sum_{l=0}^{\infty}\bigg(D_{m,\lambda}(l)+\sum_{i=0}^{l}\binom{l}{i}D_{m,\lambda}(l-i)(m)_{i,\lambda}\bigg)\frac{t^{l}}{l!} \nonumber \\
&=\ \sum_{n=0}^{\infty}\bigg(\sum_{l=0}^{n}\sum_{i=0}^{l}\binom{l}{i}\binom{n}{l}(n-l)!(-\lambda)^{n-l}(m)_{i,\lambda}D_{m,\lambda}(l-i)(1+\delta_{0,i})\bigg)\frac{t^{n}}{n!},\nonumber
\end{align*}
where $\delta_{n,k}$ is the Kronecker's symbol. \par
Therefore, by \eqref{52}, we obtain the following theorem.
\begin{theorem}
For $n\ge 0$, we have
\begin{displaymath}
D_{m,\lambda}(n+1)= \sum_{l=0}^{n}\sum_{i=0}^{l}\binom{l}{i}\binom{n}{l}(n-l)!(-\lambda)^{n-l}(m)_{i,\lambda}D_{m,\lambda}(l-i)(1+\delta_{0,i}),
\end{displaymath}
where $\delta_{n,k}$ is the Kronecker's symbol.
\end{theorem}
From \eqref{31}, we note that
\begin{equation}
m^{n}\bigg(\frac{x-1}{m}\bigg)_{n}=\sum_{k=0}^{n}V_{m,\lambda}(n,k)(x)_{k,\lambda},\quad (n\ge 0). \label{53}
\end{equation}
Thus, by \eqref{53}, we get
\begin{equation}
V_{m,\lambda}(n,0)=m^{n}\bigg(\frac{-1}{m}\bigg)_{n}=(-1)^{n}(m+1)(2m+1)\cdots((n-1)m+1).\label{53-1}
\end{equation}
From Theorem 5 and \eqref{53-1}, we have
\begin{align}
&\sum_{n=k}^{\infty}V_{m,\lambda}(n,k)\frac{1}{m^{n}}\frac{t^{n}}{n!}=\frac{1}{k!}\big(\log_{\lambda}(1+t)^{\frac{1}{m}}\big)^{k}(1+t)^{-\frac{1}{m}} \label{54}\\
&\ =\ \frac{1}{k!}\Big(\log_{\lambda}\big(\big(e^{\frac{1}{m}\log(1+t)}-1\big)+1\big)\Big)^{k}(1+t)^{-\frac{1}{m}} \nonumber \\
&\ =\ \sum_{l=k}^{\infty}S_{1,\lambda}(l,k)\frac{1}{l!}\big(e^{\frac{1}{m}\log(1+t)}-1\big)^{l}(1+t)^{-\frac{1}{m}} \nonumber \\
&\ =\ \sum_{l=k}^{\infty}S_{1,\lambda}(l,k)\sum_{j=l}^{\infty}\frac{S_{2}(j,l)}{m^{j}}\frac{1}{j!}\big(\log(1+t)\big)^{j}(1+t)^{-\frac{1}{m}} \nonumber \\
&\ =\ \sum_{j=k}^{\infty}\sum_{l=k}^{j}S_{1,\lambda}(l,k)\frac{S_{2}(j,l)}{m^{j}}\sum_{i=j}^{\infty}S_{1}(i,j)\frac{t^{i}}{i!}(1+t)^{-\frac{1}{m}}\nonumber \\
&\ =\ \sum_{i=k}^{\infty}\bigg(\sum_{j=k}^{i}\sum_{l=k}^{j}S_{1,\lambda}(l,k)\frac{S_{2}(j,l)}{m^{j}}S_{1}(i,j)\bigg)\frac{t^{i}}{i!}\sum_{p=0}^{\infty}\binom{-\frac{1}{m}}{p}t^{p}\nonumber \\
&\ =\ \sum_{i=k}^{\infty}\bigg(\sum_{j=k}^{i}\sum_{l=k}^{j}S_{1,\lambda}(l,k)\frac{S_{2}(j,l)}{m^{j}}S_{1}(i,j)\bigg)\frac{t^{i}}{i!}\sum_{p=0}^{\infty}\frac{1}{m^{p}}V_{m,\lambda}(p,0)\frac{t^{p}}{p!}\nonumber\\
&\ =\ \sum_{n=k}^{\infty}\bigg(\sum_{i=k}^{n}\sum_{j=k}^{i}\sum_{l=k}^{j}\binom{n}{i}m^{i-j}S_{1,\lambda}(l,k)S_{2}(j,l)S_{1}(i,j)V_{m,\lambda}(n-i,0)\bigg)\frac{t^{n}}{m^{n}n!}. \nonumber
\end{align}
Therefore, by comparing the coefficients on both sides of \eqref{54}, we obtain the following theorem.
\begin{theorem}
	For $n,k\ge 0$ with $n\ge k$, we have
	\begin{displaymath}
		V_{m,\lambda}(n,k)= \sum_{i=k}^{n}\sum_{j=k}^{i}\sum_{l=k}^{j}\binom{n}{i}m^{i-j}S_{1,\lambda}(l,k)S_{2}(j,l)S_{1}(i,j)V_{m,\lambda}(n-i,0).
	\end{displaymath}
\end{theorem}
By \eqref{53}, we get
\begin{align}
&\sum_{i=0}^{n}V_{m,\lambda}(n,i)(-1)^{i}\langle x\rangle_{i,\lambda}=\sum_{k=0}^{n}V_{m,\lambda}(n,k)(-x)_{k,\lambda}=m^{n}\bigg(\frac{-x-1}{m}\bigg)_{n}\label{55} \\
&\ =\ m^n \sum_{k=0}^{n}S_{1,\frac{\lambda}{m}}(n,k)\bigg(\frac{-x-1}{m}\bigg)_{k,\frac{\lambda}{m}}=\sum_{k=0}^{n}m^{n-k}S_{1,\frac{\lambda}{m}}(n,k)(-1)^k \langle x+1 \rangle_{k,\lambda} \nonumber \\
&\ =\ \sum_{k=0}^{n}m^{n-k}S_{1,\frac{\lambda}{m}}(n,k)(-1)^k\sum_{i=0}^{k}\binom{k}{i}\langle x\rangle_{i,\lambda}\langle 1\rangle_{k-i,\lambda}	\nonumber \\
&=\ \sum_{i=0}^{n}\bigg(\sum_{k=i}^{n}(-1)^{k-i}\binom{k}{i}m^{n-k}S_{1,\frac{\lambda}{m}}(n,k)\langle 1\rangle_{k-i,\lambda}\bigg)(-1)^i \langle x\rangle_{i,\lambda}\nonumber.
\end{align}
Therefore, by \eqref{55}, we obtain the following theorem.
\begin{theorem}
For $0\le i\le n$, we have
\begin{displaymath}
V_{m,\lambda}(n,i)=\sum_{k=i}^{n}(-1)^{k-i}\binom{k}{i}m^{n-k}S_{1,\frac{\lambda}{m}}(n,k)\langle 1\rangle_{k-i,\lambda}.
\end{displaymath}
\end{theorem}
From \eqref{53}, we note that
\begin{align}
m^{n}\bigg(\frac{x}{m}\bigg)_{n}\ &=\ \sum_{i=0}^{n}V_{m,\lambda}(n,i)(x+1)_{i,\lambda} \label{56} \\
&=\ \sum_{i=0}^{n}V_{m,\lambda}(n,i)\sum_{k=0}^{i}\binom{i}{k}(1)_{i-k,\lambda}(x)_{k,\lambda}\nonumber \\
&=\ \sum_{k=0}^{n}\bigg(\sum_{i=k}^{n}\binom{i}{k}V_{m,\lambda}(n,i)(1)_{i-k,\lambda}\bigg)(x)_{k,\lambda}\nonumber.	
\end{align}
On the other hand,
\begin{align}
m^{n}\bigg(\frac{x}{m}\bigg)_{n}\ &=\ m^{n}\sum_{k=0}^{n}S_{1,\lambda/m}(n,k)\bigg(\frac{x}{m}\bigg)_{k,\frac{\lambda}{m}}\label{57} \\
&=\ m^{n}\sum_{k=0}^{n}S_{1,\lambda/m}(n,k)\frac{1}{m^{k}}(x)_{k,\lambda}\nonumber \\
&=\ \sum_{k=0}^{n}m^{n-k}S_{1,\lambda/m}(n,k)(x)_{k,\lambda}\nonumber.
\end{align}
Therefore, by \eqref{56} and \eqref{57}, we obtain the inversion formula of Theorem 19.
\begin{theorem}
	For $0\le k\le n$, we have
	\begin{displaymath}
		m^{n-k}S_{1,\lambda/m}(n,k)=\sum_{i=k}^{n}\binom{n}{i}V_{m,\lambda}(n,i)(1)_{i-k,\lambda.}
	\end{displaymath}
\end{theorem}
For $0\le k\le n$, by Theorem 12, we get
\begin{align}
&W_{m+1,\lambda}(n ,k)\ =\ \frac{1}{k!(m+1)^{k}}\sum_{l=0}^{k}\binom{k}{l}(-1)^{k-l}\big(l(m+1)+1\big)_{n,\lambda}\label{58} \\
&=\ \frac{(m+1)^{n-k}}{k!}\sum_{l=0}^{k}\binom{k}{l}(-1)^{k-l}\bigg(l+\frac{1}{m+1}\bigg)_{n,\frac{\lambda}{m+1}}\nonumber \\
&=\ \frac{(m+1)^{n-k}}{k!}\sum_{l=0}^{k}\binom{k}{l}(-1)^{k-l}\bigg(l+\frac{1}{m}-\frac{1}{m(m+1)}\bigg)_{n,\frac{\lambda}{m+1}}\nonumber \\
&=\ (m+1)^{n}\sum_{l=0}^{k}\binom{k}{l}\frac{(-1)^{k-l}}{(m+1)^{k}{k!}}\sum_{j=0}^{n}\binom{n}{j}\bigg(l+\frac{1}{m}\bigg)_{n-j,\frac{\lambda}{m+1}}\bigg(\frac{-1}{m(m+1)}\bigg)_{j,\frac{\lambda}{m+1}} \nonumber\\
&=\ (m+1)^{n}\sum_{l=0}^{k}\binom{k}{l}\frac{(-1)^{k-l}}{(m+1)^{k}k!}\sum_{j=0}^{n}\binom{n}{j}\frac{(lm+1)_{n-j,\frac{m}{m+1}\lambda}}{m^{n-j}}\bigg(\frac{-1}{m(m+1)}\bigg)_{j,\frac{\lambda}{m+1}} \nonumber \\
&=\ (m+1)^{n}\sum_{j=0}^{n}\binom{n}{j}\frac{m^{k}}{(m+1)^{k}m^{n-j}}\bigg(-\frac{1}{m(m+1)}\bigg)_{j,\frac{\lambda}{m+1}}\nonumber\\
&\quad\quad\times \frac{1}{k!m^{k}}\sum_{l=0}^{k}\binom{k}{l}(-1)^{k-l}(lm+1)_{n-j,\frac{m}{m+1}\lambda}\nonumber \\
&=\ (m+1)^{n}\sum_{j=0}^{n}\binom{n}{j}\frac{m^{k}}{(m+1)^{k}m^{n-j}}\frac{(-1)^{j}}{(m(m+1))^{j}}\langle 1\rangle_{j,m\lambda}W_{m,\frac{m}{m+1}\lambda}(n-j,k)\nonumber \\
&=\ \frac{1}{(m+1)^{k}m^{n-k}}\sum_{j=0}^{n}\binom{n}{j}(-1)^{n-j}(m+1)^{j}\langle 1\rangle_{n-j,m\lambda}W_{m,\frac{m}{m+1}\lambda}(j,k)\nonumber.	
\end{align}
Therefore, by \eqref{58}, we obtain the following theorem.
\begin{theorem}
For $0\le k\le n$, we have
\begin{displaymath}
W_{m+1,\lambda}(n,k)=\frac{1}{(m+1)^{k}m^{n-k}}\sum_{j=0}^{n}\binom{n}{j}(-1)^{n-j}(m+1)^{j}\langle 1\rangle_{n-j,m\lambda}W_{m,\frac{m}{m+1}\lambda}(j,k).
\end{displaymath}
\end{theorem}
From \eqref{26} and Theorem 21, we can derive the following equation:
\begin{align}
& D_{m+1,\lambda}(n,x)\ =\ \sum_{k=0}^{n}W_{m+1,\lambda}(n,k)x^{k}\label{59} \\
&=\ \sum_{k=0}^{n}\bigg(\frac{1}{(m+1)^{k}m^{n-k}}\sum_{j=0}^{n}\binom{n}{j}(-1)^{n-j}(m+1)^{j}\langle 1\rangle_{n-j,m\lambda}W_{m,\lambda}(j,k)\bigg)x^{k}\nonumber 	\\
&=\ \frac{1}{m^{n}}\sum_{j=0}^{n}\binom{n}{j}(-1)^{n-j}(m+1)^{j}\langle 1\rangle_{n-j,m\lambda}\sum_{k=0}^{n}W_{m,\lambda}(j,k)\bigg(\frac{m}{m+1}x\bigg)^{k}\nonumber \\
&=\ \frac{1}{m^{n}}\sum_{j=0}^{n}\binom{n}{j}(-1)^{n-j}(m+1)^{j}\langle 1\rangle_{n-j,m\lambda}D_{m,\lambda}\bigg(j,\frac{m}{m+1}x\bigg). \nonumber
\end{align}
Here $W_{m,\lambda}(j,k)=0$, if $k>j$. \par
Therefore, we obtain the following corollary.
\begin{corollary}
For $n\ge 0$, we have
\begin{displaymath}
D_{m+1,\lambda}(n,x)= \frac{1}{m^{n}}\sum_{j=0}^{n}\binom{n}{j}(-1)^{n-j}(m+1)^{j}\langle 1\rangle_{n-j,m\lambda}D_{m,\lambda}\bigg(j,\frac{m}{m+1}x\bigg).
\end{displaymath}
\end{corollary}
\noindent \emph{Remark.} Note that
\begin{align*}
	\mathcal{F}_{m+1,\lambda}(n,x)&=\sum_{k=0}^{n}k!W_{m+1,\lambda}(n,x)x^{k}\\
	&=\ \sum_{k=0}^{n}\bigg(\frac{k!}{(m+1)^{k}m^{n-k}}\sum_{j=0}^{n}\binom{n}{j}(-1)^{n-j}(m+1)^{j}\langle 1\rangle_{n-j,m\lambda}W_{m,\lambda}(j,k)\bigg)x^{k}\\
	&=\ \frac{1}{m^{n}}\sum_{j=0}^{n}\binom{n}{j}(-1)^{n-j}(m+1)^{j}\langle 1\rangle_{n-j,m\lambda}\sum_{k=0}^{n}k!W_{m,\lambda}(j,k)\bigg(\frac{m}{m+1}x\bigg)^{k} \\
	&=\ \frac{1}{m^{n}}\sum_{j=0}^{n}\binom{n}{j}(-1)^{n-j}(m+1)^{j}\langle 1\rangle_{n-j,m\lambda}\mathcal{F}_{m,\lambda}\bigg(j,\frac{m}{m+1}x\bigg).
\end{align*}
Now, we observe from Theorem 13 that
\begin{align}
D_{m,\lambda}(n,x)\ &=\ \sum_{k=0}^{n}W_{m,\lambda}(n,k)x^{k}\ =\ \sum_{k=0}^{n}\sum_{i=k}^{n}\binom{n}{i}m^{i-k}(1)_{n-i,\lambda}S_{2,\lambda/m}(i,k)x^{k}\label{60} \\
&=\ \sum_{i=0}^{n}\binom{n}{i}m^{i}(1)_{n-i,\lambda}\sum_{k=0}^{i}S_{2,\lambda/m}(i,k)\bigg(\frac{x}{m}\bigg)^{k}\nonumber \\
&=\ \sum_{i=0}^{n}\binom{n}{i}m^{i}(1)_{n-i,\lambda}\mathrm{Bel}_{i,\lambda/m}\bigg(\frac{x}{m}\bigg).\nonumber
\end{align}
Therefore, by \eqref{60}, we obtain the following theorem.
\begin{theorem}
	For $n\ge 0$, we have
	\begin{displaymath}
		D_{m,\lambda}(n,x)=\sum_{i=0}^{n}\binom{n}{i}m^{i}(1)_{n-i,\lambda}\mathrm{Bel}_{i,\lambda/m}\bigg(\frac{x}{m}\bigg).
	\end{displaymath}
\end{theorem}
 From the definition of $\lambda$-falling factorial sequences, we have
 \begin{align}
 (t)_{n,\lambda}\ &=\ (t-1+1)_{n,\lambda}\ =\ \sum_{k=0}^{n}\binom{n}{k}(1)_{n-k,\lambda}(t-1)_{k,\lambda}\label{61} \\
 &=\ \sum_{k=0}^{n}\binom{n}{k}(1)_{n-k,\lambda}\sum_{j=0}^{k}\binom{k}{j}(-1)_{k-j,\lambda}(t)_{j,\lambda}\nonumber \\
 &=\ \sum_{k=0}^{n}\binom{n}{k}(1)_{n-k,\lambda}\sum_{j=0}^{k}\binom{k}{j}(-1)^{k-j}\langle 1\rangle_{k-j,\lambda}(t)_{j,\lambda}\nonumber \\
 &=\ \sum_{j=0}^{n}\bigg(\sum_{k=j}^{n}\binom{n}{k}(1)_{n-k,\lambda}(-1)^{k-j}\langle 1\rangle_{k-j,\lambda}\binom{k}{j}\bigg)(t)_{j,\lambda}.\nonumber
 \end{align}
In the same way, we also have
\begin{align}
(t)_{n,\lambda}\ &=\ (t+1-1)_{n,\lambda}\ =\ \sum_{k=0}^{n}\binom{n}{k}(-1)^{n-k}\langle 1\rangle_{n-k,\lambda}(t+1)_{k,\lambda}\label{61-1}\\
&=\ \sum_{k=0}^{n}\binom{n}{k}(-1)^{n-k}\langle 1\rangle_{n-k,\lambda}\sum_{j=0}^{k}\binom{k}{j}(1)_{k-j,\lambda}(t)_{j,\lambda}\nonumber\\
&=\ \sum_{j=0}^{n}\bigg(\sum_{k=j}^{n}\binom{n}{k}\binom{k}{j}(1)_{k-j,\lambda}\langle 1\rangle_{n-k,\lambda}(-1)^{n-k}\bigg)(t)_{j,\lambda}. \nonumber
\end{align}
Thus, from \eqref{61} and \eqref{61-1} we have the following lemma.
\begin{lemma}
For $n\ge 0$, we have
\begin{displaymath}
\sum_{k=j}^{n}\binom{n}{k}\binom{k}{j}(-1)^{k-j}(1)_{n-k,\lambda}\langle 1\rangle_{k-j,\lambda}=\delta_{n,j}=\sum_{k=j}^{n}\binom{n}{k}\binom{k}{j}(-1)^{n-k}\langle 1\rangle_{n-k,\lambda}(1)_{k-j,\lambda}.
\end{displaymath}
\end{lemma}
Suppose that $\displaystyle a_{n,\lambda}=\sum_{k=0}^{n}\binom{n}{k}(1)_{n-k,\lambda}b_{k,\lambda}\displaystyle$. Then, by Lemma 24, we have
\begin{align}
&\sum_{k=0}^{n}\binom{n}{k}(-1)^{n-k}\langle 1\rangle_{n-k,\lambda}a_{k,\lambda}=\sum_{k=0}^{n}\binom{n}{k}(-1)^{n-k}\langle 1\rangle_{n-k,\lambda}\sum_{j=0}^{k}\binom{k}{j}(1)_{k-j,\lambda}b_{j,\lambda}\label{62} 	\\
&=\ \sum_{j=0}^{n}\bigg(\sum_{k=j}^{n}\binom{n}{k}\binom{k}{j}(-1)^{n-k}(1)_{k-j,\lambda}\langle 1\rangle_{n-k,\lambda}\bigg)b_{j,\lambda} \nonumber \\
&=\ b_{n,\lambda}.\nonumber
\end{align}
Conversely, we assume that $\displaystyle b_{n,\lambda}=\sum_{k=0}^{n}\binom{n}{k}(-1)^{n-k}\langle 1\rangle_{n-k,\lambda}a_{k,\lambda}$. Then, by Lemma 24, we note that
\begin{align}
\sum_{k=0}^{n}\binom{n}{k}(1)_{n-k,\lambda}b_{k,\lambda}\ &=\ \sum_{k=0}^{n}\binom{n}{k}(1)_{n-k,\lambda}\sum_{j=0}^{k}\binom{k}{j}(-1)^{k-j}\langle 1\rangle_{k-j,\lambda}a_{j,\lambda} \label{63} 	\\
&=\ \sum_{j=0}^{n}\bigg(\sum_{k=j}^{n}\binom{n}{k}\binom{k}{j}(-1)^{k-j}\langle 1\rangle_{k-j,\lambda}(1)_{n-k,\lambda}\bigg)a_{j,\lambda}\nonumber \\
&=\ a_{n,\lambda}. \nonumber
\end{align}
Therefore, by \eqref{62} and \eqref{63}, we obtain the following theorem.
\begin{theorem}
For $n\ge 0$, we have
\begin{displaymath}
a_{n,\lambda}=\sum_{k=0}^{n}\binom{n}{k}(1)_{n-k,\lambda}b_{k,\lambda}\ \Longleftrightarrow\ b_{n,\lambda}=\sum_{k=0}^{n}\binom{n}{k}(-1)^{n-k}\langle 1\rangle_{n-k,\lambda}a_{k,\lambda}.
\end{displaymath}
\end{theorem}
From Theorem 23 and Theorem 25, we obtain the following theorem.
\begin{theorem}
For $n\ge 0$, we have
\begin{displaymath}
m^{n}\mathrm{Bel}_{n,\lambda/m}\bigg(\frac{x}{m}\bigg)=\sum_{k=0}^{n}\binom{n}{k}(-1)^{n-k}\langle 1\rangle_{n-k,\lambda}D_{m,\lambda}(k,x).
\end{displaymath}
\end{theorem}

\section{Degenerate $r$-Whitney numbers of Dowling lattices}
For $r\in\mathbb{N}$, we consider the degenerate $r$-Whitney numbers of the first kind $V_{m,\lambda}^{(r)}(n,k)$ and those of the second $W_{m,\lambda}^{(r)}(n,k)$, respectively defined by
\begin{equation}
m^{n}(x)_{n}=\sum_{k=0}^{n}V_{m,\lambda}^{(r)}(n,k)(mx+r)_{k,\lambda}, \label{64}
\end{equation}
and
\begin{equation}
(mx+r)_{n,\lambda}=\sum_{k=0}^{n}m^{k}W_{m,\lambda}^{(r)}(n,k)(x)_{k},\quad(n\ge 0). \label{65}
\end{equation}
We remark that $r$-Whitney numbers and their applications were studied by several authors (see\ [3,4,20]). Note that $\displaystyle\lim_{\lambda\rightarrow 0}V_{m,\lambda}^{(r)}(n,k)=V_{m}^{(r)}(n,k) \displaystyle$ and  $\displaystyle\lim_{\lambda\rightarrow 0}W_{m,\lambda}^{(r)}(n,k)=W_{m}^{(r)}(n,k) \displaystyle$ are respectively the $r$-Whitney numbers of the first kind and those of the second kind.
From \eqref{64}, we note that
\begin{align}
(1+mt)^{x}\ &=\ \sum_{n=0}^{\infty}(x)_{n}m^{n}\frac{t^{n}}{n!}=\sum_{n=0}^{\infty}\bigg(\sum_{k=0}^{n}V_{m,\lambda}^{(r)}(n,k)(mx+r)_{k,\lambda}\bigg)\frac{t^{n}}{n!}\label{66} \\
&=\ \sum_{k=0}^{\infty}\bigg(\sum_{n=k}^{\infty}V_{m,\lambda}^{(r)}(n,k)\frac{t^{n}}{n!}\bigg)(mx+r)_{k,\lambda}.\nonumber	
\end{align}
On the other hand,
\begin{align}
(1+mt)^{x}\ &=\ (1+mt)^{x+\frac{r}{m}}(1+mt)^{-\frac{r}{m}}\ =\ (1+mt)^{\frac{mx+r}{m}}(1+mt)^{-\frac{r}{m}} \label{67} \\
&=\ e_{\lambda}^{mx+r}\big(\log_{\lambda}(1+mt)^{\frac{1}{m}}\big)(1+mt)^{-\frac{r}{m}}\nonumber \\
&=\ \sum_{k=0}^{\infty}\bigg(\frac{1}{k!}\Big(\log_{\lambda}(1+mt)^{\frac{1}{m}}\Big)^{k}(1+mt)^{-\frac{r}{m}}\bigg)(mx+r)_{k,\lambda}.\nonumber
\end{align}
By \eqref{66} and \eqref{67}, we get the generating function of the degenerate $r$-Whitney numbers of first kind given by
\begin{equation}
\frac{1}{k!}\Big(\log_{\lambda}e_{m}(t)\Big)^{k}e_{m}^{-r}(t)=\sum_{n=k}^{\infty}V_{m,\lambda}^{(r)}(n,k)\frac{t^{n}}{n!},\label{68}
\end{equation}
where $k$ is a nonnegative integer. \par
From \eqref{65}, we have
\begin{align}
e_{\lambda}^{mx+r}(t)\ &=\ \sum_{n=0}^{\infty}(mx+r)_{n,\lambda}\frac{t^{n}}{n!}\ =\ \sum_{n=0}^{\infty}\bigg(\sum_{k=0}^{n}m^{k}W_{m,\lambda}^{(r)}(n,k)(x)_{k}\bigg)\frac{t^{n}}{n!}\label{69} \\
&=\ \sum_{k=0}^{\infty}\bigg(\sum_{n=k}^{\infty}m^{k}W_{m,\lambda}^{(r)}(n,k)\frac{t^{n}}{n!}\bigg)(x)_{k}.\nonumber
\end{align}
On the other hand,
\begin{equation}
e_{\lambda}^{mx+r}=e_{\lambda}^{r}(t)\big(e_{\lambda}^{m}(t)-1+1\big)^{x}=\sum_{k=0}^{\infty}\bigg(\frac{1}{k!}\big(e_{\lambda}^{m}(t)-1\big)^{k}e_{\lambda}^{r}(t)\bigg)(x)_{k}.\label{70}	
\end{equation}
By \eqref{69} and \eqref{70}, we get
\begin{equation}
\frac{1}{k!}\bigg(\frac{e_{\lambda}^{m}(t)-1}{m}\bigg)^{k}e_{\lambda}^{r}(t)=\sum_{n=k}^{\infty}W_{m,\lambda}^{(r)}(n,k)\frac{t^{n}}{n!},\quad(k\ge 0).\label{71}
\end{equation}
It is known that the unsigned degenerate $r$-Stirling numbers of the first kind are defined by
\begin{equation}
\langle x+r\rangle_{n}=\sum_{k=0}^{n}{n+r \brack k+r}_{r,\lambda}\langle x\rangle_{k,\lambda},\quad(n\ge 0),\quad(\mathrm{see}\ [18]), \label{72}
\end{equation}
where $\langle x\rangle_{0}=1,\ \langle x\rangle_{n}=x(x+1)(x+2)\cdots(x+n-1),\ (n\ge 1)$. \par
From \eqref{72}, we have
\begin{equation}
(1-t)^{-r}\frac{1}{k!}\big(-\log_{\lambda}(1-t)\big)^{k}=\sum_{n=k}^{\infty}{n+r \brack k+r}_{r,\lambda}\frac{t^{n}}{n!},\label{73}	
\end{equation}
where $k$ is a nonnegative integer (see [18]). \par
By \eqref{73} and \eqref{68}, we easily get
\begin{equation}
\sum_{n=k}^{\infty}V_{1,\lambda}^{(r)}(n,k)\frac{t^{n}}{n!}=(1+t)^{-r}\frac{1}{k!}\big(\log_{\lambda}(1+t)\big)^{k}=\sum_{n=k}^{\infty}(-1)^{n-k}{n+r \brack k+r}_{r,\lambda}\frac{t^{n}}{n!}. \label{74}	
\end{equation}
Comparing the coefficients on both sides of \eqref{74}, we have
\begin{equation}
V_{1,\lambda}^{(r)}(n,k)= (-1)^{n-k}{n+r \brack k+r}_{r,\lambda},\quad V_{1,\lambda}^{(0)}=S_{1,\lambda}(n,k), \quad V_{m,\lambda}^{(1)}(n,k)=V_{m,\lambda}(n,k). \label{75}
\end{equation}
In \cite{18}, the degenerate $r$-Stirling numbers of the second kind are defined by 
\begin{equation}
(x+r)_{n,\lambda}=\sum_{k=0}^{n}{n+r \brace k+r}_{r,\lambda}(x)_{k},\quad(n\ge 0).\label{76}
\end{equation}
From \eqref{76}, we can easily derive the generating function of the degenerate $r$-Stirling number of the second kind given by
\begin{equation}
\frac{1}{k!}\big(e_{\lambda}(t)-1\big)^{k}e_{\lambda}^{r}(t)=\sum_{n=k}^{\infty}{n+r \brace k+r}_{r,\lambda}\frac{t^{n}}{n!}, \label{77}	
\end{equation}
where $k$ is a nonnegative integer (see \cite{15}). \par
From \eqref{71} and \eqref{77}, we have
\begin{displaymath}
W_{1,\lambda}^{(r)}(n,k)={n+r \brace k+r}_{r,\lambda},\quad W_{1,\lambda}^{(0)}(n,k)=S_{2,\lambda}(n,k),\quad W_{m,\lambda}^{(1)}(n,k)=W_{m,\lambda}(n,k).
\end{displaymath}

\section{Further remarks}
From \eqref{19}, we have
\begin{align}
\sum_{n=0}^{\infty}\beta_{n,\lambda}^{(k)}\frac{t^{n}}{n!}\ &=\ \bigg(\frac{t}{e_{\lambda}(t)-1}\bigg)^{k}\label{80} \\
&=\ \frac{k!}{(e_{\lambda}(t)-1)^{k}}\frac{1}{k!}\Big(\log_{\lambda}\big((e_{\lambda}(t)-1\big)+1\Big)^{k} \nonumber \\
&=\  \frac{k!}{(e_{\lambda}(t)-1)	}\sum_{l=k}^{\infty}\frac{(e_{\lambda}(t)-1)^{l}}{l!}S_{1,\lambda}(l,k)\nonumber \\
&=\ \frac{k!}{(e_{\lambda}(t)-1)^{k}}\sum_{l=0}^{\infty}S_{1,\lambda}(l+k,k)\frac{(e_{\lambda}(t)-1)^{l+k}}{(l+k)!}\nonumber \\
&=\ \sum_{l=0}^{\infty}S_{1,\lambda}(l+k,k)\frac{l!k!}{(l+k)!}\frac{1}{l!}\big(e_{\lambda}(t)-1\big)^{l} \nonumber \\
&=\ \sum_{l=0}^{\infty}S_{1,\lambda}(l+k,k)\frac{1}{\binom{l+k}{l}}\sum_{n=l}^{\infty}S_{2,\lambda}(n,l)\frac{t^{n}}{n!} \nonumber\\
&=\ \sum_{n=0}^{\infty}\bigg(\sum_{l=0}^{n}\binom{l+k}{l}^{-1}S_{1,\lambda}(l+k,k)S_{2,\lambda}(n,l)\bigg)\frac{t^{n}}{n!}. \nonumber
\end{align}
From \eqref{20}, we obtain
\begin{align}
\sum_{n=0}^{\infty}\mathcal{E}_{n,\lambda}^{(\alpha)}\frac{t^{n}}{n!}\  &=\ \bigg(\frac{2}{e_{\lambda}(t)+1}\bigg)^{\alpha}\ =\ \bigg(\frac{e_{\lambda}(t)-1}{2}+1\bigg)^{-\alpha}\label{81}	\\
&=\ \sum_{l=0}^{\infty}(-1)^{l}\binom{\alpha+l-1}{l}\bigg(\frac{e_{\lambda}(t)-1}{2}\bigg)^{l}\nonumber \\
& \ =\ \sum_{l=0}^{\infty}\bigg(-\frac{1}{2}\bigg)^{l}\binom{\alpha+l+1
}{l}l!\frac{1}{l!}\big(e_{\lambda}(t)-1\big)^{l}\nonumber
\end{align}
\begin{align*}
&=\ \sum_{l=0}^{\infty}\bigg(-\frac{1}{2}\bigg)^{l}\binom{\alpha+l+1}{l}l!\sum_{n=l}^{\infty}S_{2,\lambda}(n,l)\frac{t^{n}}{n!}\nonumber \\
&=\ \sum_{n=0}^{\infty}\bigg(\sum_{l=0}^{n}\bigg(-\frac{1}{2}\bigg)^{l}\binom{\alpha+l-1}{l}l!S_{2,\lambda}(n,l)\bigg)\frac{t^{n}}{n!}. \nonumber
\end{align*}
We summarize our results in \eqref{80} and \eqref{81} as a theorem.

\begin{theorem}
For $n\ge 0$, we have
\begin{align*}
\beta_{n,\lambda}^{(k)}=\sum_{l=0}^{n}\binom{l+k}{k}^{-1}S_{1,\lambda}(l+k,k)S_{2,\lambda}(n,l), \\
\mathcal{E}_{n,\lambda}^{(\alpha)}=\sum_{l=0}^{n}\bigg(-\frac{1}{2}\bigg)^{l}\binom{\alpha+l-1}{l}l!S_{2,\lambda}(n,l).
\end{align*}
\end{theorem}

\end{document}